\documentclass[12pt]{amsart}  
\usepackage[utf8]{inputenc}
\usepackage[english]{babel}
\usepackage{graphicx}
\usepackage{siunitx}
\usepackage{graphicx}
\usepackage{verbatim}
\usepackage{relsize}
\usepackage{amsthm}
\usepackage{mathtools}
\usepackage{mathrsfs}
\usepackage{blindtext}
\usepackage{float}
\usepackage{geometry} 
\usepackage{amsmath}
\usepackage{amssymb}
\usepackage{latexsym}
\usepackage{cancel}
\usepackage{tikz}
\usepackage{lipsum}
\usepackage{tikz-cd}
\usepackage{wasysym}
\usepackage[mathscr]{eucal}

\usepackage{layout}
\usepackage{enumerate}
\usepackage{amsfonts}
\usepackage[geometry]{ifsym}
\usetikzlibrary{tikzmark,shapes.misc}
\title{Function spaces on Corson-like compacta}

\author[K.\ Zakrzewski]{Krzysztof Zakrzewski}
\address{Institute of Mathematics\\
University of Warsaw\\ Banacha 2\newline 02--097 Warszawa\\
Poland\\}
\email{Krzysztof.Zakrzewski@mimuw.edu.pl}
\thanks{The author was partially supported by the NCN
(National Science Centre, Poland) research grant no.\ 2020/37/B/ST1/02613.
}
\theoremstyle{plain}
\newtheorem{Theorem}{Theorem}[section]
\newtheorem{Lemma}[Theorem]{Lemma}
\newtheorem{Definition}[Theorem]{Definition}
\newtheorem{prop}[Theorem]{Proposition}

\newtheorem{Example}[Theorem]{Example}
\numberwithin{equation}{section}
\DeclareMathOperator{\Limm}{Lim}
\newcommand{\Lim}[1]{\raisebox{0.5ex}{\scalebox{0.8}{$\displaystyle \lim_{#1}\;$}}}
\newcommand\norm[1]{\left\lVert#1\right\rVert}

\newsavebox{\overlongequation}

\begin{document}

\begin{abstract}
For an index set $\Gamma$ and a cardinal number $\kappa$ the $\Sigma_{\kappa}$-product of real lines $\Sigma_{\kappa}(\mathbb{R}^{\Gamma})$ consist of all elements of $\mathbb{R}^{\Gamma}$ with $<\kappa$ nonzero coordinates. A compact space is $\kappa$-Corson if it can be embedded into $\Sigma_{\kappa}(\mathbb{R}^{\Gamma})$ for some $\Gamma$. We also consider a class of compact spaces wider than the class of $\omega$-Corson compact spaces, investigated by Nakhmanson and Yakovlev as well as Marciszewski, Plebanek and Zakrzewski called $NY$ compact spaces. For a Tychonoff space $X$, let $C_{p}(X)$ be the space of real continuous functions on the space $X$, endowed with the pointwise convergence topology.

We present here a characterisation of $\kappa$-Corson compact spaces $K$ for regular, uncountable cardinal numbers $\kappa$ in terms of function spaces $C_{p}(K)$, extending a theorem of Bell and Marciszewski and a theorem of Pol.
 We also prove that classes of $NY$ compact spaces and $\omega$-Corson compact spaces $K$  are preserved by linear homeomorphisms of function spaces $C_{p}(K)$. 
\end{abstract}
\maketitle
\section{Introduction}
For an infinite cardinal number $\kappa$, a compact space $K$ is called $\kappa$-Corson compact if it is homeomorphic to a subset of the $\Sigma_{\kappa}$-product of real lines for some set $\Gamma$, that is
$$\Sigma_{\kappa}(\Gamma)=\{x\in \mathbb{R}^{\Gamma}: |\{\gamma\in\Gamma:x(\gamma)\neq 0\}|<\kappa \}.$$
For $\kappa=\omega_{1}$ the class of $\kappa$-Corson compact spaces is the well known class of Corson compact spaces, see \cite{N}, \cite{AMN}. Classes of $\kappa$-Corson compact spaces, for arbitrary $\kappa$, were investigated in \cite{K} and \cite{BM}.

The classes of $\omega$-Corson compact spaces and $NY$ compact spaces were studied in a paper by Marciszewski, Plebanek and Zakrzewski, see \cite{MPZ}, which encouraged the author to investigate function spaces $C_{p}(K)$ for such $K$. Here is the summary of main results:
\begin{enumerate}[(A)]

\item In Section \ref{secdimension} we prove that, for $\sigma$-compact spaces $X$ and $Y$, where $X$ is strongly countable dimensional, if function spaces $C_{p}(X)$ and $C_{p}(Y)$ are linearly homeomorphic, then $Y$ is strongly countable dimensional as well (Theorem \ref{sil kappa wymiar zach przy lin homeo}). We also investigate the case when $X$ is finite dimensional (Theorem \ref{sk wymiar}).  
\item In Section \ref{secNY} we show that $NY$ compactness is invariant under linear homeomorphisms of function spaces. In fact we prove several slightly more general theorems (Theorems \ref{NY zachowuja sie przy lin homeo}, \ref{1 tw o zan eb w NY}, \ref{2 tw o zan w NY}). We also give an example of an Eberlein compact space $K$ and a \mbox{$\omega$-Corson} compact space $L$ such that $C_{p}(K)$ is linearly homeomorphic to a linear subspace of $C_{p}(L)$ but $K$ is not $NY$ compact (Example \ref{przykl Eberleina}).

\item In Section \ref{char kappa-cors sekcja} we prove a characterisation  (Theorem \ref{char kappa - Corson}) of $\kappa$-Corson compact spaces $K$ in terms of function spaces $C_{p}(K)$ for regular uncountable cardinal numbers $\kappa$,  which extends results of Pol \cite{P} and Bell and Marciszewski \cite{BM}. Additionally, we show that, for compact spaces $K$ and $L$, if $K$ is \mbox{$\omega$-Corson} compact, and function spaces $C_{p}(K)$ and $C_{p}(L)$ are linearly homeomorphic, then the space $L$ is $\omega$-Corson compact as well (Theorem \ref{omega corson zach przy lin homeo}). 
\end{enumerate}
\section{Preliminaries}
\subsection{Terminology and notation}
Let $\mathbb{N}=\{1,2,\dots\}$ be the set of natural numbers.
All topological spaces under consideration are assumed to be completely regular (Tychonoff).
For a topological space $X$, let $C_{p}(X)$ denote the space of real continuous functions on $X$ endowed with the pointwise convergence topology.
As usual, the topological weight of a space $X$ is denoted by $w(X)$.

A family $\mathcal{U}$ of subsets of a space $X$
is called \emph{$T_0$-separating if}, for every pair of distinct points $x,y$ of
$X$, there is $U\in\mathcal{U}$ containing exactly one of the
points $x,y$.

Given a family $\mathcal{U}$ of subsets of a space $X$ we define $ord(x,\mathcal{U}) = |\{U\in \mathcal{U}: x\in U\}|$ for $x\in X$ and $ord(\mathcal{U})=sup\{ord(x,\mathcal{U}):x\in X\}$. We say that $\mathcal{U}$ is \emph{point-finite} if $ord(x,\mathcal{U}) < \omega$ for all $x\in X$.

A space $X$ is called \textit{metacompact} if every open cover $\mathcal{U}$ of $X$ has a point-finite, open refinement $\mathcal{V}$. It is \textit{$\sigma$-metacompact} if every open cover $\mathcal{U}$ of $X$ has an open refinement which is a union of countably many point-finite families.
A space is said to be hereditarily metacompact ($\sigma$-metacompact) if its every subspace is metacompact ($\sigma$-metacompact).

For a space $X$, sets of the form $f^{-1}((0,1])$, where $f:X \to [0,1]$ is a continuous function, are called the \textit{cozero sets}.

By (finite) dimension of a space $X$ we will always mean the covering dimension $dim(X)$. Recall that a space $X$ is called \emph{strongly countable dimensional} if $X$ is a countable union of closed subspaces of finite covering dimension.

For a Banach space $X$, let $B_{X}$ denote the closed unit ball of $X$.

For a space $X$ let $A(X)$ denote the one-point compactification of the space $X$.

\subsection{$\Sigma_{\kappa}$-products}\label{subsec_sigma_prod}
In the sequel, $\kappa$ always stands for an infinite cardinal number.

For $x$ in the product space $\mathbb{R}^\Gamma$, let
$supp(x)=\{\gamma\in\Gamma: x(\gamma)\neq 0\}$.
 The $\Sigma_{\kappa}$-product is defined as
$$\Sigma_\kappa(\mathbb{R}^\Gamma)=\{x\in\mathbb{R}^\Gamma: |supp(x)|<\kappa\}.$$

The space $\Sigma_\kappa ([0,1]^\Gamma)$ is defined analogously.

We will introduce more general notation for the case when $\kappa=\omega$.
Let $\{X_\gamma: \gamma\in \Gamma\}$ be a family of nonempty topological spaces, and let $a_\gamma\in X_{\gamma}$. The $\sigma$-product of the family
$\{(X_\gamma,a_\gamma): \gamma\in \Gamma\}$ is defined as
$$\sigma(X_\gamma,a_\gamma,\Gamma) = \Big\{(x_\gamma)_{\gamma\in \Gamma}\in \prod_{\gamma\in \Gamma} X_\gamma: |\{\gamma\in\Gamma: x_\gamma \ne a_\gamma\}| < \omega\Big\}\,.$$

For the special cases when $X_\gamma = I$ and $a_\gamma = 0$, for every $\gamma\in \Gamma$ or $X_\gamma = I^{\omega}$ and $a_\gamma = 0$, for every $\gamma\in \Gamma$  we denote
$\sigma(X_\gamma,a_\gamma,\Gamma)$ by $\sigma(I,\Gamma)$ or $\sigma(I^{\omega},\Gamma)$ respectively.
\subsection{$\kappa$-Corson compact spaces and $NY$ compact spaces}

For an infinite cardinal number $\kappa$, a compact space is called $\kappa$-Corson compact if it is homeomorphic to a subset of $\Sigma_{\kappa}(\mathbb{R}^{\Gamma})$ for some set $\Gamma$. Equivalently, instead of $\mathbb{R}^{\Gamma}$ one can consider $I^{\Gamma}$. Every compact space is $\kappa$-Corson compact for $\kappa>w(K)$. Traditionally, instead of $\omega_{1}$-Corson compact space we will write Corson compact space.

A compact space is said to be $NY$ compact if it is homeomorphic to a subset of $\sigma(X_{\gamma},a_{\gamma},\Gamma)$ for some compact, metrizable $X_{\gamma}$'s, $a_{\gamma}\in X_{\gamma}$ and some set $\Gamma$. Clearly, every $\omega$-Corson compact is $NY$ compact. It is easy to observe that a compact space is $NY$ compact if and only if it is homeomorphic to a subspace of $\sigma(I^{\omega},\Gamma)$ for some set $\Gamma$ (see \cite[Proposition 3.2]{MPZ}). The class of $NY$ compact spaces is invariant under taking continuous images.

Recall that a space is Eberlein compact if it is homeomorphic to a compact subset of a Banach space endowed with the weak topology. Equivalently, by  the celebrated Amir-Lindenstrauss theorem, a compact space is Eberlein compact if and only if it is homeomorphic to a subset of  
$$c_{0}(\Gamma)=\{x\in\mathbb{R}^{\Gamma}:\,|\{\gamma\in\Gamma:|x(\gamma)|<\epsilon\}|<\omega\text{  for every }\epsilon>0\}.$$
Since $\sigma(I^{\omega},\Gamma)\sim\sigma(\prod_{n\in\mathbb{N}}[0,\frac{1}{n}],0,\Gamma)\subseteq c_{0}(\Gamma\times\mathbb{N})$, every $NY$ compact space is Eberlein compact.

\subsection{Linear topological spaces and dual operators}
\label{operators}
For a real linear topological space $X$, let $X^{*}$ denote its dual space consisting of continuous linear functionals $\phi:X\xrightarrow[]{}\mathbb{R}$.
Recall that, for a continuous linear operator $T:X\xrightarrow[]{}Y$ between two linear topological vector spaces, the dual operator $T^{*}:Y^{*}\xrightarrow[]{}X^{*}$ is given by the formula $T^{*}(\phi)=\phi\circ T$.

For two locally convex linear topological spaces $X$ and $Y$, a continuous operator \\$T:X\xrightarrow[]{}Y$ has a dense image if and only if its dual operator is injective (\cite[Corollary 21/Chapter 26]{S}).

Let us note here that, from \cite[Theorem 7.3]{Sc} and  \cite[Corollary 21/Chapter 26]{S} it follows that, for two locally convex linear topological spaces $X$ and $Y$, a continuous linear operator $T:X\xrightarrow[]{}Y$ is a weakly open (onto its image) injection if and only if the dual operator is a surjection.

For a topological space $X$, the space $C_{p}(X)$ is a locally convex linear topological space.
\section{Behaviour of dimension under transformations of function spaces}
\label{secdimension}
In 1980 Pavlovski proved that the covering dimension of metrizable, compact spaces is preserved under linear homeomorphisms of function spaces endowed with the pointwise convergence topology (\cite{Pv}). Later in 1982 Pestov showed that this holds for arbitrary Tychonoff spaces (see \cite{Pe}). In 1990 Arhangel'skii asked if it is true that $dim\,Y\leq dim\, X$ when there exists a continuous, linear surjection of $C_{p}(X)$ onto $C_{p}(Y)$. In 1997 Leiderman, Morris and Pestov showed that the answer to this question is in the negative. Namely, they proved that for every  metrizable, finite dimensional compact space $X$, there exists a continuous, linear surjection from $C_{p}(I)$ to $C_{p}(X)$  (see \cite{LMP}).   

In 1997 Leiderman, Levin and Pestov proved that, for every finite dimensional, metrizable, compact space $Y$, there exists a 2-dimensional, metrizable, compact space $X$ such that there is a continuous, open, linear surjection of $C_{p}(X)$ onto $C_{p}(Y)$. They also proved that, for every $n\in\mathbb{N}$, there exist an $n$-dimensional, compact, metrizable space $Y$ and 1-dimensional, compact, metrizable space $X$ such that there is a continuous, linear, open surjection from $C_{p}(X)$ to $C_{p}(Y)$.  Lastly, they gave an example of a continuous, linear surjection from $C_{p}([0,1])$ to  $C_{p}([0,1])$ which is not open.  For details see \cite{LLP}.

In the same year Levin and Sternfeld noticed that one can construct a continuous surjection of $C_{p}(I)$ onto $C_{p}(X)$ where $X$ is an infinite dimensional space (see \cite[Remark 4.6]{LMP}).  

In 2011 Levin proved a strengthening of results from \cite{LMP} and \cite{LLP}, namely he showed that for every metrizable, finite dimensional, compact space $X$, the space $C_{p}(X)$ is a continuous, linear, open image of $C_{p}([0,1])$ (see \cite{L}).

For a compact space $K$, define $I(K)=\bigcup\{U\subseteq K:\text{U is open and finite dimensional}\}$. Let $K^{(0)}=K$, $K^{(\alpha+1)}=K^{(\alpha)}\setminus I(K)$ for every ordinal number $\alpha$ and $K^{(\lambda)}=\bigcap_{\alpha<\lambda}K^{(\alpha)}$ for a limit ordinal $\lambda$. The fd-height of a compact space $K$ is the minimal ordinal number $\alpha$ such that $K^{(\alpha)}=\emptyset$. A compact space is strongly countable dimensional if and only if its fd-height is a countable, ordinal number.

In 2017 Gartside and Feng proved that if $K$ is a compact, metrizable space and has finite fd-height, then there exists a  continuous, linear surjection from $C_{p}(I)$ to $C_{p}(K)$. They also showed that if $X$ is submetrizable and has $k_{\omega}$ sequence of compact sets of finite fd-height, then there exists a continuous, linear surjection from $C_{p}(\mathbb{R})$ to $C_{p}(X)$. Finally, they showed that if there is a continuous, linear surjection of $C_{p}(\mathbb{R})$ onto $C_{p}(X)$, then $X$ is submetrizable, has $k_{\omega}$ sequence and is strongly countable dimensional (for details and definitions see \cite{GF}). 

Also in 2017 Kawamura and Leiderman proved that for compact spaces $K$ and $L$ if there exists a continuous linear surjection $T:C_{p}(K)\xrightarrow[]{}C_{p}(L)$, then if $K$ is zero-dimensional, then $L$ is zero-dimensional as well (see \cite{KL}).

In 2019 Górak, Krupski and Marciszewski proved that a space $X$ is compact, metrizable and strongly countable dimensional if and only if for every $\epsilon>0$ there exists a $(1+\epsilon)$-good, uniformly continuous surjection $u_{\epsilon}:C_{p}(I)\xrightarrow[]{}C_{p}(X)$ such that $\norm{u_{\epsilon}(f)}_{\infty}\leq\norm{f}_{\infty}$ for every $f\in C_{p}(I)$. They also proved that if $X$ is compact, metrizable, and there exists a uniformly continuous surjection from $C_{p}(X)$ to $C_{p}(Y)$ which is $c$-good for some $c>0$, then if $X$ is zero-dimensional, then $Y$ is zero-dimensional, and if $X$ is strongly countable dimensional, then so is $Y$ (for details and definitions see \cite{GKM}).

The following two lemmas are \cite[Proposition 6.7.2]{M} and \cite[Proposition 6.7.4]{M}.

    For a space $X$, let $L_{p}(X)$ denote the dual space of $C_{p}(X)$. 

\begin{Lemma}
\label{X zan sie w Lp}
Every space $X$ embeds as a closed subset in $L_{p}(X)$.
\end{Lemma}
\begin{Lemma}\label{funkcjonały na Cp}
For a space $X$, every $\phi\in L_{p}(X)$ is of the form 
$$\phi(f)=\sum_{k=1}^{n}a_{k}f(x_{k})$$
where $a_{k}\in\mathbb{R}$, $x_{k}\in X$ and $n\in\mathbb{N}$.
\end{Lemma}
\begin{Definition}
    For a space $X$ and $\phi\in L_{p}(X)$, let $supp(\phi)=\{x_{1},\dots,x_{n}\}\subseteq X$ if\\
    $\phi(f)=\sum_{k=1}^{n}a_{k}f(x_{k})$ for every $f\in C_{p}(X)$ and some nonzero $a_{1},\dots, a_{n}\in\mathbb{R}$.
\end{Definition}

\begin{Lemma} Let $X$ be a space, and let $n\in\mathbb{N}$. The set 
$$A_{n}=\{\phi\in L_{p}(X):|supp(\phi)|\leq n\}$$
is closed in $L_{p}(X)$ .
\label{funkcjonały o nośnikach <=n są zbiorem domkniętym}
\end{Lemma}
\begin{proof}
Pick any $\phi\in L_{p}(X)\setminus A$, and let $supp(\phi)=\{x_{1},\dots,x_{m}\}$ where $m\geq n+1$. There exists pairwise disjoint open sets $U_{i}\subseteq X$ such that $U_{i}\cap supp(\phi)=\{x_{i}\}$ for every $i\in\{1,\dots,m\}$. Since the space $X$ is Tichonoff, there exists $f_{i}\in C_{p}(X)$ such that $f_{i}\restriction_{X\setminus U_{i}}\equiv 0$ and $f_{i}(x_{i})=1$ for every $i\in\{1,\dots,m\}$. The set
$$V=\{\phi\in L_{p}(X):\phi(f_{i})\neq 0,\text{ for }i\in\{1,\dots,m\}\}$$
is an open neighbourhood of $\phi$ disjoint from $A$.
\end{proof}
The following theorem and lemma are \cite[Theorem 3.3.7]{E} and \cite[Lemma 3.1.6]{E}.
\begin{Theorem}[Theorem on dimension rising mappings]
If $f:X\xrightarrow[]{}Y$ is a closed mapping from a normal space $X$ onto a normal space $Y$ and there exists $k\in\mathbb{N}$ such that $|f^{-1}(y)|\leq k$ for every $y\in Y$, then $dim\, Y\leq dim\, X + (k-1)$.
\end{Theorem}
\begin{Lemma}\label{uogólnienie wymiaru sumy}
If a normal space $X$ can be represented as a union of a sequence $K_{1},K_{2},\dots$ of subspaces such that for every $i\in\mathbb{N}$ and every $F\subseteq K_{i}$ closed in $X$ we have $dim\,F\leq n$ and the union $\bigcup_{j\leq i}K_{j}$  is closed for every $i\in\mathbb{N}$, then $dim\, X\leq n$.
\end{Lemma}
\begin{Lemma}
\label{wymiar ilorazu}
Let $X$ be a finite dimensional, normal topological space, and let $A\subseteq X$ be a closed subset, then $dim\, X/A\leq dim\, X$.
\end{Lemma}
\begin{proof}
Let $\pi:X\xrightarrow[]{}X/A$ be the quotient map.
Let $\mathcal{U}$ be a finite open cover of $X/A$. Consider
$$\mathcal{V}=\{\pi^{-1}(U):U\in\mathcal{U}\},$$
it is a open cover of $X$. There exists $V_{0}\in\mathcal{V}$ such that $A\subseteq V_{0}$. Let
$$\mathcal{V}^{'}=\{V_{0}\}\cup\{V\setminus A:V\in\mathcal{V}\setminus\{V_{0}\}\},$$ it is an open cover of $X$. By \cite[Theorem 3.2.1]{E} there exists a shrinking $\mathcal{W}$ of $\mathcal{V}^{'}$ such that $ord(\mathcal{W})\leq dim\,X+1$ . Notice that there is precisely one $W_{0}\in \mathcal{W}$ such that $W_{0}\cap A\neq\emptyset$, and so for this $W_{0}$, we have $A\subseteq W_{0}$. The family
$$\mathcal{U}^{'}=\{\pi(W):W\in\mathcal{W}\}$$
is an open cover of $X/A$ inscribed in $U$ and $ord(\mathcal{U}^{'})=ord(\mathcal{W})\leq dim X+1$. Consequently, we have $dim\, X/A\leq dim\,X.$
\end{proof}

\begin{Definition}
    A topological space $X$ is called strongly countable dimensional if it can be represented as a union of countably many closed, finite dimensional subspaces.
\end{Definition}
\begin{Definition}
 For  $\kappa\geq\omega$, we say that a topological space $X$ is strongly $\kappa$-dimensional if it can be represented as a union of $\kappa$ many closed, finite dimensional subspaces.
\end{Definition}
\begin{Definition}
For $\kappa\geq\omega$, we say that a topological space $X$ is $\kappa$-compact if it is a union of its $\kappa$ many compact subspaces.
\end{Definition}
A theorem of Gartside and Feng shows the thesis of the theorem below for the case when $T$ is a surjection, $X=\mathbb{R}$ and $\kappa=\omega$. Recall that by a result of Górak, Krupski and Marciszewski for a metrizable, compact $X$ and $\kappa=\omega$, it is enough to assume that $T$ is a uniformly continuous surjection which is $c$-good for some $c>0$.
\begin{Theorem}
Let $\kappa\in Card$, $\kappa\geq \omega$. Assume that $X$ and $Y$ are $\kappa$-compact Tychonoff spaces. Assume further there exists a continuous, linear transformation  $T:C_{p}(X)\xrightarrow[]{}C_{p}(Y)$ such that the image of $T$ is dense in $C_{p}(Y)$.
 If $X$ is strongly \mbox{$\kappa$-dimensional,} then $Y$ is strongly $\kappa$-dimensional as well. \label{sil kappa wymiar zach przy lin homeo}
\end{Theorem}
\begin{proof}
It is easy to notice that $X=\bigcup_{\alpha\in\kappa}X_{\alpha}$ where $X_{\alpha}$'s are compact, finite dimensional, and the family $\{X_{\alpha}:\alpha\in \kappa\}$ is closed under finite unions.
Let $$s_{n,m}^{\alpha}:X^{n}_{\alpha}\times ([-m,-1/m]\cup[1/m,m])^{n}\xrightarrow[]{}L_{p}(X)$$
be given by the formula
$$s_{n,m}^{\alpha}(x_{1},\dots,x_{n},a_{1},\dots,a_{n})(f)=a_{1}f(x_{1})+\dots+a_{n}f(x_{n}).$$
Since $\{X_{\alpha}:\alpha\in\kappa\}$ is closed under finite unions, we have
$$\{\mathbf{0}\}\cup\bigcup\{im(s^{\alpha}_{n,m}):\alpha\in\kappa,n,m\in\mathbb{N}\}=L_{p}(X).$$
The function $s_{n,m}^{\alpha}$ is closed as a continuous function from a compact space.
Define 
$$S_{n}=\{\phi\in L_{p}(X):|supp(\phi)|<n\}$$
and
$$Y_{\alpha,n}=\{x\in X_{\alpha}^{n}:x_{k}=x_{l} \text{ for some }k\neq l\},$$
then
$$Y_{\alpha,n}\times([-m,-1/m]\cup[1/m,m])^{n}=(s_{n,m}^{\alpha})^{-1}(S_{n}).$$
Clearly, the space $Y_{\alpha,n}$ is a closed subspace of $X^{n}_{\alpha}$, and by Lemma \ref{funkcjonały o nośnikach <=n są zbiorem domkniętym}, the space $S_{n}$ is closed in $L_{p}(X)$.
Put 
$$Z_{\alpha,n,m}=X^{n}_{\alpha}\times([-m,-1/m]\cup[1/m,m])^{n}/(Y_{\alpha,n}\times ([-m,-1/m]\cup[1/m,m])^{n}).$$
By Lemma \ref{wymiar ilorazu}, we have $dim\,(Z_{\alpha,n,m})\leq n\cdot dim\,X_{\alpha}+n$, moreover $Z_{\alpha,n,m}$ is compact and Hausdorff. 
Let $p_{n,m}^{\alpha}:Z_{\alpha,n,m}\xrightarrow[]{}L_{p}(X)/S_{n}$ be given by the formula
$$p_{n,m}^{\alpha}([x])=\begin{cases}
[s_{n,m}^{\alpha}(x)] & \text{if } x\notin Y_{\alpha,n}\times([-m,-1/m]\cup[1/m,m])^{n}\\
[\mathbf{0}] & \text{if } x\in Y_{\alpha,n}\times([-m,-1/m]\cup[1/m,m])^{n} 
\end{cases}.$$
It is a continuous function from a compact, Hausdorff space into a Hausdorff space, so it is closed. For every $\phi\in L_{p}(X)$, we have $|(p_{n,m}^{\alpha})^{-1}([\phi])|\leq n!$, since $|(p^{\alpha}_{n,m})^{-1}([\mathbf{0}])|=1$ and $|(p^{\alpha}_{n,m})^{-1}([\phi])|=n!$ for $\phi\notin S_{n}$. By the theorem on dimension rising mappings, we have 
$$dim\,im (p^{\alpha}_{n,m})\leq dim\,Z_{\alpha,n,m}+n!-1\leq n\cdot dim\,X_{\alpha}+n+n!-1,$$
that is
$$dim\,im(s^{\alpha}_{n,m})/(im(s^{\alpha}_{n,m})\cap S_{n})\leq n\cdot dim\,X_{\alpha}+n+n!-1.$$

We will prove by induction on $n$ that, for every $n,m\in\mathbb{N}$ and $\alpha\in\kappa$, we have
$$dim\,im(s^{\alpha}_{n,m})\leq n\cdot dim\,X_{\alpha}+n+n!-1.$$
Firstly, as $s^{\alpha}_{1,m}$ is a homeomorphic embedding for every $m\in\mathbb{N}$ and $\alpha\in\kappa$, and $S_{1}=\{\mathbf{0}\}$, we get that 
$$dim\,im(s_{1,m}^{\alpha})/(im(s_{1,m}^{\alpha})\cap S_{1})=dim\,im(s_{1,m}^{\alpha})=dim\,X_{\alpha}\times[-m,m]\leq dim\,X_{\alpha}+1$$
for every $m\in\mathbb{N}$ and $\alpha\in\kappa$. This shows the thesis for $n=1$.
Assume that for every $k<n,m\in\mathbb{N}$ and $\alpha\in\kappa$ 
$$dim\,im(s^{\alpha}_{k,m})\leq k\cdot dim\,X_{\alpha}+k+k!-1=D(\alpha,k).$$
By the inductive assumption, the space $\bigcup_{m\in\mathbb{N}}\bigcup_{k\leq n-1}im(s^{\alpha}_{k,m})$ is a union of countably many compact spaces of dimension at most $D(\alpha,n-1)$, so 
$$dim\,\bigcup\{im(s^{\alpha}_{k,m}):m\in\mathbb{N}, k\leq n-1\}\leq D(\alpha,n-1).$$
Obviously,
$$im(s^{\alpha}_{n,m})=(im(s^{\alpha}_{n,m})\cap S_{n})\cup(im(s^{\alpha}_{n,m})\setminus S_{n}).$$

Firstly, notice that
$im(s^{\alpha}_{n,m})\cap S_{n}$is a closed subset of $\bigcup_{m\in\mathbb{N}}\bigcup_{k\leq n-1}im(s^{\alpha}_{k,m})$, and therefore $im(s^{\alpha}_{n,m})\cap S_{n}$ is a normal space of dimension at most $D(\alpha,n-1)$.

Secondly, let $F\subseteq im(s^{\alpha}_{n,m})\setminus S_{n}$ be closed in $im(s^{\alpha}_{n,m})$, and let
$$\pi^{\alpha}_{n,m}:im(s^{\alpha}_{n,m})\xrightarrow[]{}im(s^{\alpha}_{n,m})/(im(s^{\alpha}_{n,m})\cap S_{n})$$
be the quotient map. The set $\pi_{n,m}^{\alpha}(F)\sim F$ is closed in $im(s^{\alpha}_{n,m})/(im(s^{\alpha}_{n,m})\cap S_{n})$, and therefore
$$dim\, F=dim\,\pi^{\alpha}_{n,m}(F)\leq dim\,im(s^{\alpha}_{n,m})/(im(s^{\alpha}_{n,m})\cap S_{n})\leq n\cdot dim\,X_{\alpha}+n+n!-1=D(\alpha,n)$$
By Lemma \ref{uogólnienie wymiaru sumy}, we have 
$$dim\,im(s^{\alpha}_{n,m})\leq max(D(\alpha,n-1),D(\alpha,n))=D(\alpha,n).$$
As 
$$L_{p}(X)=\{\mathbf{0}\}\cup\bigcup\{im(s^{\alpha}_{n,m}):\alpha\in\kappa, n,m\in\mathbb{N}\},$$
the space $L_{p}(X)$ is strongly $\kappa$-dimensional. By \cite[Corollary 21/Chapter 26]{S}, there exists a continuous, linear injection  $T^{*}:L_{p}(Y)\xrightarrow[]{}L_{p}(X)$.  By Lemma \ref{X zan sie w Lp}, the space $Y$ is homeomorphic to a closed subspace $Y'$ of $L_{p}(Y)$. Let $Y'=\bigcup_{\alpha\in\kappa}Y_{\alpha}$ where $Y_{\alpha}$'s are compact spaces. For $\alpha\in\kappa$, the restriction $T^{*}\restriction_{Y_{\alpha}}$ is a homeomorphic embedding, so $Y_{\alpha}$ is strongly $\kappa$-dimensional. Since $Y\sim Y'=\bigcup_{\alpha\in\kappa}Y_{\alpha}$, the space $Y$ is strongly $\kappa$-dimensional as well. 
\end{proof}
In the theorem below, to have a bound on dimension of $Y$, the assumption that support is bounded is necessary because  by a theorem of Levin, for every metrizable, finite dimensional, compact space $K$, the space $C_{p}(K)$ is a continuous, linear, open image of $C_{p}([0,1])$. For the standard definition of $supp_{T}(y)$ see \cite{M}. Equivalently, $supp_{T}(y)=supp(T^{*}(e(y)))$.
\begin{Theorem}\label{sk wymiar}
Let $X$ and $Y$ be $\sigma$-compact spaces such that $X$ is finite dimensional. If there exists a continuous, linear transformation $T:C_{p}(X)\xrightarrow[]{}C_{p}(Y)$ such that $T(C_{p}(X))$ is dense in $C_{p}(Y)$, and  $|supp_{T}(y)|\leq p$ for some $p\in\mathbb{N}$ and every $y\in Y$, then $$dim(Y)\leq p\cdot dim\,X+p+p!-1 .$$
\end{Theorem}
\begin{proof}
Using definitions from the proof of Theorem \ref{sil kappa wymiar zach przy lin homeo} for $\kappa=\omega$, we have 
$$X=\bigcup_{k\in\mathbb{N}}X_{k}$$
where $X_{k}$'s are compact, and
$$L_{p}(X)=\bigcup\{im(s^{k}_{n,m}):k,n,m\in\mathbb{N}\}\cup\{\mathbf{0}\}.$$ Let $Y'\subseteq L_{p}(Y)$, $Y'\sim Y$ be the subspace of evaluations (Lemma \ref{X zan sie w Lp}).
Since $|supp_{T}(y)|\leq p$ for every $y\in Y$, there exists a continuous, linear injection 
$$T^{*}\restriction_{Y'}:Y'\xrightarrow[]{}\{\phi\in L_{p}(X):|supp(\phi)|\leq p\}=\bigcup\{im(s_{n,m}^{k}):k,m\in\mathbb{N},n\leq p\}\cup\{\mathbf{0}\}.$$
Let $Y'=\bigcup_{l\in\mathbb{N}}Y_{l}$ for some compact spaces $Y_{l}$. For every $l\in\mathbb{N}$, the restriction $T^{*}\restriction_{Y_{l}}$ is a homeomorphic embedding.
We have 
$$dim\,im(s^{k}_{n,m})\leq n\cdot dim\,X+n+n!-1$$ for every $k,n,m\in\mathbb{N}$, so
$$dim\,Y_{l}\leq dim\,\bigcup_{k,m\in\mathbb{N},n\leq p}im(s_{n,m}^{k})\cup\{\mathbf{0}\}\leq p\cdot dim\,X+p+p!-1$$
for every $l\in\mathbb{N}$ by Lemma \ref{uogólnienie wymiaru sumy}. Again by Lemma \ref{uogólnienie wymiaru sumy}, we can conclude that \[dim\,Y\leq p\cdot dim\,X+p+p!-1.\]
\end{proof}
\section{Function spaces on $NY$ compact spaces }
\label{secNY}
The main result of this section is the invariance of the class of $NY$ compact spaces $K$ under linear homeomorphisms of function spaces $C_{p}(K)$. In fact we prove several more general theorems (Theorems \ref{NY zachowuja sie przy lin homeo}, \ref{1 tw o zan eb w NY}, \ref{2 tw o zan w NY}). Our reasoning relies on the condition $(iii)$ in the following characterisation proved in \cite{MPZ}. The equivalence of the conditions $(i)$ and $(ii)$ was proved in \cite{NY}.
\begin{Theorem}
   \label{charact_ec_sigma_m}
For a compact space $K$, the following
conditions are equivalent:
\begin{enumerate}[(i)]
\item $K$ is $NY$ compact;
\item there exists a $T_0$-separating family 
$\mathcal{U} = \bigcup\{\mathcal{U}_\gamma: \gamma\in \Gamma\}$ 
consisting of cozero subsets of $K$, where each $\mathcal{U}_\gamma$ is countable and the family $\{\bigcup\mathcal{U}_\gamma: \gamma\in \Gamma\}$ is point-finite;
\item $K$ is hereditarily metacompact, and each nonempty subspace $A$ of $K$ 
contains a nonempty, relatively open subspace $U$ of countable weight.
\end{enumerate}
\end{Theorem}
\begin{Lemma}
Let $K$ be a compact space such that $K=\bigcup_{n\in\mathbb{N}}K_{n}$ where
$\{K_{n}:n\in\mathbb{N}\}$ is a sequence of compact spaces with the property that  every nonempty subspace $A$ of $K_{n}$ contains a nonempty, relatively open subspace $U$ of countable weight. Then the space $K$ has the same property.
\label{pzel suma NY rel otw lem}
\end{Lemma}
\begin{proof}
Assume there exists $A\subseteq K$ such that every relatively open $U\subseteq A$ has uncountable weight, then $\overline{A}$ has the same property. Indeed, assume there is a 
nonempty, relatively open $U\subseteq \overline{A}$ of countable weight, 
then $U\cap A\neq \emptyset$, and therefore $U\cap A\subseteq A$ is a relatively open, nonempty subset of $A$ with $w(U\cap A)\leq \omega$, contradiction. 

Without loss of generality, we can assume that $A$ is closed and therefore compact. Then it has the same property as $K$, so we can assume that $A=K$. As $K=\bigcup_{n\in\mathbb{N}}K_{n}$, by the Baire category theorem, there is $K_{n}$ with nonempty interior. By the assumption, there exists a nonempty, open $V\subseteq int_{K}(K_{n})$ with $w(V)\leq \omega$. Then $V$ is a nonempty, open subset of $K$ of countable weight.
\end{proof}

\begin{Lemma}
    Let $K$ be an $NY$ compact space, and let $L$ be a compact space. If there exists a continuous, linear  transformation $T:C_{p}(K)\xrightarrow[]{} C_{p}(L)$ such that $T(C_{p}(K))$ is dense in $C_{p}(L)$, then $L$ is a countable union of $NY$ compact spaces.
    \label{przel suma NY zach sie przy lin homeo}
\end{Lemma}
\begin{proof}
For $n,m\in\mathbb{N}$, let 
$$s_{n,m}:K^{n}\times [-m,m]^{n}\xrightarrow[]{}L_{p}(K)$$
be given by the formula
$$s_{n,m}((x_{1},\dots, x_{n}), (t_{1},\dots, t_{n}))=\sum_{k=1}^{n}t_{k}e(x_{k})$$
where $e(x)(f)=f(x)$ for every $f\in C_{p}(K)$ and $x\in K$. Functions $s_{n,m}$ are continuous.  By Lemma \ref{funkcjonały na Cp}, we have
$$\bigcup\{im(s_{n,m}):n,m\in\mathbb{N}\}=L_{p}(K).$$

A finite product of $NY$ compact spaces is again $NY$ compact, so  $K^{n}\times[-m,m]^{n}$ is $NY$ compact for $n,m\in\mathbb{N}$. As $NY$ compact spaces are preserved under continuous images, the space $L_{p}(K)$ is a union of countably many $NY$ compact spaces. By \cite[Corollary 21/Chapter 26]{S}, there exists a continuous, linear injection  $T^{*}:L_{p}(L)\xrightarrow[]{}L_{p}(K)$, and by Lemma \ref{X zan sie w Lp}, the space $L$ is homeomorphic to a subspace  $L'$ of $L_{p}(L)$. Since $L$ is a compact space, the transformation $T^{*}\restriction_{L'}$ is a homeomorphic embedding. We obtain that $L$ is a union of countably many $NY$ compact spaces.
\end{proof}
Using Lemmas \ref{pzel suma NY rel otw lem} and \ref{przel suma NY zach sie przy lin homeo} we obtain
\begin{prop}
Let $K$ be an $NY$ compact space, and let $L$ be a compact space. Assume there exists a continuous, linear transformation $T:C_{p}(K)\xrightarrow[]{} C_{p}(L)$ such that $T(C_{p}(K))$ is dense in $C_{p}(L)$. Then every nonempty subspace  $A$ of $L$ contains a nonempty, relatively open subspace of countable weight.
\label{X-NY i Cp(X) lin homeo Cp(Y) to Y ma wl podprzestrzeni}
\end{prop}
The above proposition shows that the space $L$ follows the second part of condition $iii)$ from Theorem \ref{charact_ec_sigma_m}. Now it suffices to show that $L$ is hereditarily metacompact.

The following lemma is \cite[Lemma 5.3.5]{E}.
\begin{Lemma}
For every open cover $\{U_{s}:s\in S\}$ of a metacompact space, there exists a point-finite open cover $\{V_{s}:s\in S\}$ such that $V_{s}\subseteq U_{s}$ for every $s\in S$.\label{En}
\end{Lemma}
The following lemma is a generalisation of \cite[Lemma 5.3.6]{E} for closed transformations $f$, and is proved in an analogous way.
\begin{Lemma}
    Let $f:X\xrightarrow[]{}Y$ be a continuous map from a metacompact space $X$ onto space $Y$ such that the image of every closed subset $F$ of $X$ is an $F_{\sigma}$ set in $Y$. Then for every open cover  $\mathcal{U}$ of $Y$ which is a union of countably many point-finite families, there exists a point-finite open refinement $\mathcal{V}$. 
    \label{uogólnienie lematu z EN}
\end{Lemma}
\begin{proof}
Let $\mathcal{U}=\bigcup_{n=1}^{\infty}\mathcal{U}_{n}$ where $\mathcal{U}_{n}$'s are point-finite families. By Lemma \ref{En}, there exists a point-finite open cover $\{G_{n}:n\in\mathbb{N}\}$ of the space $X$ such that $G_{n}\subseteq f^{-1}(U_{n})$ where $U_{n}=\bigcup \mathcal{U}_{n}$. The set $E_{n}=X\setminus \bigcup_{k\geq n} G_{k}$ is closed for $n\in \mathbb{N}$. One can easily see that $E_{1}\subseteq E_{2}\subseteq \dots$. Since the family $\{G_{n}:n\in\mathbb{N}\}$ is point finite, the family $\{E_{n}:n\in\mathbb{N}\}$ is a cover of $X$. Moreover, we have
$$f(E_{n})\subseteq f\Big(\bigcup_{k<n}G_{k}\Big)= \bigcup_{k<n}f(G_{k})\subseteq \bigcup_{k<n}U_{k},\text{  for } n\in\mathbb{N}.$$
For $n\in\mathbb{N}$, let $f(E_{n})=\bigcup_{k\in\mathbb{N}}E_{n}^{k}$ where $E_{n}^{k}$ are closed. Without loss of generality we can assume that $E_{n}^{k}\subseteq E_{n}^{k+1}$ and $E_{n}^{k}\subseteq E_{n+1}^{k}$ for $k,n\in\mathbb{N}$. For $n\in\mathbb{N}$, let
$$\mathcal{V}_{n}=\{U\setminus E_{n}^{n}:U\in\mathcal{U}_{n}\}.$$
It suffices to show that $\mathcal{V}=\bigcup_{n\in\mathbb{N}}\mathcal{V}_{n}$ is a point-finite cover of $Y$. We shall prove that $\mathcal{V}$ is a cover of $Y$.

Pick $y\in Y$, and let $n(y)$ be the smallest number such that $y\in U_{n(y)}$. Since we have $f(E_{n(y)})\subseteq \bigcup_{n<n(y)}U_{n}$, there exists $U\in\mathcal{U}_{n(y)}$ such that $y\in U\setminus f(E_{n(y)})$, then $y\in U\setminus E_{n(y)}^{k}$ for every $k\in\mathbb{N}$, in particular $y\in U\setminus E_{n(y)}^{n(y)}\in \mathcal{V}_{n(y)}$. 

It remains to show that $\mathcal{V}$ is point-finite. Choose any $y\in Y$. As $Y=\bigcup_{n\in\mathbb{N}}f(E_{n})$, the point $y$ belongs to a set $f(E_{n})$ for some $n\in\mathbb{N}$. Consequently, $y\in E_{n}^{k}$ for some $k\in\mathbb{N}$, and if $m=max(n,k)$, then $y\in E_{m}^{m}$. By the definition of $\mathcal{V}_{n},$ we  have $y\notin \bigcup_{n>m}\mathcal{V}_{n}$. Since the families $\mathcal{U}_{n}$ are point-finite, the point $y$ belongs to only finitely many elements of $\mathcal{V}$.
\end{proof}
\begin{Lemma}
    Every $\sigma$-metacompact space which is a union of its countably many closed metacompact subspaces is metacompact.
    \label{sigma-metazw + przel sum domk meta = meta}
\end{Lemma}
\begin{proof}
Assume $X$ is $\sigma$-metacompact, and $X=\bigcup_{n\in\mathbb{N}}X_{n}$ where $X_{n}\subseteq X$  are closed and metacompact. Let $\mathcal{U}$ be an open cover of $X$. By $\sigma$-metacompactness, there exists a $\sigma$-point finite, open cover $\mathcal{V}\prec \mathcal{U}$.

Consider $Y=\dot{\bigcup}_{n\in\mathbb{N}}X_{n}$ endowed with the topology of disjoint union, and let $f:Y\xrightarrow[]{}X$ be such that $f\restriction_{X_{n}}$ is the identity embedding of $X_{n}$ into $X$. Clearly, the space $Y$ is metacompact, and the transformation $f$ is a continuous surjection with the property that image of every closed set in $Y$ is an $F_{\sigma}$ set in $X$. By Lemma \ref{uogólnienie lematu z EN}, there exists a point-finite open cover $\mathcal{W}\prec\mathcal{V}\prec\mathcal{U}$. The space $X$ is metacompact.
\end{proof}
The following theorem is \cite[Theorem 7.1]{A}.
\begin{Theorem}\label{char eb}
For a compact space $K$, the following conditions are equivalent:
\begin{enumerate}[(i)]
\item $K$ is Eberlein compact;
\item $C_{p}(K)$ has a $\sigma$-compact dense subspace.
\end{enumerate}
\end{Theorem}
The following theorem can be found in \cite{G} (here $\Delta = \{(x,x): x\in K\}$ is the diagonal). The implication $(i)\implies(ii)$ was proven by Yakovlev (\cite{Ya}).
\begin{Theorem}[Gruenhage]\label{charact_Gruenhage}
	For a compact space $K$, the following
	conditions are equivalent:
\begin{enumerate}[(i)]
\item  $K$ is Eberlein compact;
\item $K^2$ is hereditarily $\sigma$-metacompact;
\item $K^2\setminus\Delta$ is  $\sigma$-metacompact.
\end{enumerate}
\end{Theorem}
\begin{Lemma}
\label{her met eb}
An Eberlein compact space which is a union of its countably many closed, hereditarily metacompact subspaces is hereditarily metacompact.
\end{Lemma}
\begin{proof}
Let $K$ be an Eberelein compact space, and let $K=\bigcup_{n\in\mathbb{N}}K_{n}$ where $K_{n}$'s are hereditarily metacompact, compact spaces. Let $A\subseteq K$. By Theorem \ref{charact_Gruenhage}, the space $A$ is $\sigma$-metacompact. The spaces $A\cap K_{n}$ are closed in $A$ and metacompact. By Lemma \ref{sigma-metazw + przel sum domk meta = meta}, the space $A$ is metacompact. This shows that $K$ is hereditarily metacompact.
\end{proof}

\begin{Theorem}
Let $K$ be an $NY$ compact space, and let $L$ be a compact space. If there exists a continuous, linear transformation $T:C_{p}(K)\xrightarrow[]{}C_{p}(L)$ such that $T(C_{p}(K))$ is dense in $C_{p}(L)$, then $L$ is $NY$ compact as well. \label{NY zachowuja sie przy lin homeo}
\end{Theorem}
\begin{proof}
By Theorem \ref{charact_ec_sigma_m} and Proposition \ref{X-NY i Cp(X) lin homeo Cp(Y) to Y ma wl podprzestrzeni}, it suffices to show that $L$ is hereditarily metacompact. By Lemma \ref{przel suma NY zach sie przy lin homeo}, we have $L=\bigcup_{n}L_{n}$ where $L_{n}$ is $NY$ compact for every $n\in\mathbb{N}$. By Theorem \ref{char eb}, the space $C_{p}(K)$ has a $\sigma$-compact dense subspace $D$, then $T(D)$ is a $\sigma$-compact dense subspace of $C_{p}(L)$, and again by Theorem \ref{char eb}, the space $L$ is Eberlein compact. By Lemma \ref{her met eb}, the space $L$ is hereditarily metacompact.
\end{proof}
\begin{Lemma}
Let $K$ and $L$ be compact spaces. Assume there exists a closed linear homeomorphic embedding $T:C_{p}(K)\xrightarrow[]{}C_{p}(L)$. If $L$ is Eberlein compact, then $K$ is Eberlein compact as well.
\label{lem do zan eb}
\end{Lemma}
\begin{proof}
Let $C(K),C(L)$ denote the Banach spaces of real function spaces on $K$ and $L$ respectively. Since the topology on $C(L)$ is stronger than the topology on $C_{p}(L)$, the image $T(C(K))$ is closed in $C(L)$ and therefore is a Banach space. By the closed graph theorem, transformation $T:C(K)\xrightarrow[]{}C(L)$ is also continuous. By the open mapping theorem, it is a homeomorphic embedding. Consequently, the dual transformation $T^{*}:C(L)^{*}\xrightarrow[]{}C(K)^{*}$ is a continuous, linear surjection. By the open mapping theorem, the function $T^{*}$ is open, so there exists $c>0$ such that $c\cdot B_{C(K)^{*}}\subseteq T^{*}(B_{C(L)^{*}})$. 

Since $L$ is Eberlein compact, the space $(B_{C(L)^{*}},\omega^{*})$ is also Eberlein compact.   Obviously $T^{*}(B_{C(L)^{*}})$ is an Eberlein compactum as a continuous image of an Eberlein compact space, then $B_{C(K)^{*}}$ is also Eberlein compact. Finally, the space $K$ is homeomorphic to a subspace of $B_{C(K)^{*}}$, so it is Eberlein compact as well.
\end{proof}
\begin{Theorem}
Let $K$ be an Eberlein compact space, and let $L$ be an $NY$ compact space. If there exists a linear homeomorphic embedding $T:C_{p}(K)\xrightarrow[]{}C_{p}(L)$, then $K$ is $NY$ compact as well.\label{1 tw o zan eb w NY}
\end{Theorem}
\begin{proof}
 From the proof of Lemma \ref{przel suma NY zach sie przy lin homeo}, we know that if $L$ is $NY$ compact, then $L_{p}(L)$ is a union of countably many $NY$ compact spaces. By \cite[Theorem 7.3]{Sc} and \cite[Corollary 21/Chapter 26]{S}, there exists a continuous surjection $T^{*}:L_{p}(L)\xrightarrow[]{}L_{p}(K)$, so $L_{p}(K)$ is a union of countably many $NY$ compact spaces as well. By Lemma \ref{X zan sie w Lp}, the space $K$ is a countable union of $NY$ compact spaces and therefore closed hereditarily metacompact spaces. By Lemma \ref{her met eb}, the space $K$ is hereditarily metacompact. By Lemma \ref{pzel suma NY rel otw lem}, every subspace of $K$ has an open set of countable weight. The space $K$ is \mbox{$NY$ compact} by Theorem \ref{charact_ec_sigma_m}.
\end{proof}
By Lemma \ref{lem do zan eb} and Theorem \ref{1 tw o zan eb w NY}, we obtain
\begin{Theorem}
Let $K$ and $L$ be compact spaces. Assume there exists a closed, linear, homeomorphic embedding $T:C_{p}(K)\xrightarrow[]{}C_{p}(L)$. If $L$ is $NY$ compact, then $K$ is $NY$ compact as well.\label{2 tw o zan w NY}
\end{Theorem}
The following example shows that assumptions in Theorems $\ref{1 tw o zan eb w NY}$, $\ref{2 tw o zan w NY}$ are necessary.
\begin{Example}
There exist a compact space $K$ and a $\omega$-Corson compact space $L$ such that there is a linear, homeomorphic embedding $T:C_{p}(K)\xrightarrow[]{}C_{p}(L)$, but $K$ is not Eberlein compact.
\label{przykl Eberleina}
\end{Example}
\begin{proof}
Let $Z=2^{\omega}\,\cup\,2^{<\omega}$, where
points from $2^{<\omega}$ are isolated in $Z$, and basic 
neighbourhoods of a point $x\in 2^{\omega}$ are of the 
form 
$$\{x\}\cup\{x\restriction_{n}:n\geq k\}, \text{ for some } k\in\omega.$$
Let $K=2^{\omega}\cup 2^{<\omega}\cup\{\infty\}$ be the one-point compactification of $Z$. Notice that $2^{<\omega}$ is a countable, dense subset of $K$. Clearly, the space $K$ is not metrizable as it has an uncountable discrete set $2^{\omega}$. This shows that $K$ is not Eberlein compact because every separable Eberlein compactum is metrizable. 

Let $L=(2^{\omega}\cup\{\infty\})\dot\cup(\omega+1)$, where $2^{\omega}\cup\{\infty\}$ is the one-point compactification of $2^{\omega}$ endowed with the discrete topology and $\omega+1$ is endowed with the order topology. Clearly, we have $L\sim A(\mathfrak{c})\,\dot\cup\, A(\omega)$, and therefore $L$ is $\omega$-Corson compact, since it is a disjoint union of two $\omega$-Corson compact spaces. 

Let $\sigma:2^{<\omega}\xrightarrow[]{}\omega$ be a bijection, and let $T:C_{p}(K)\xrightarrow[]{}C_{p}(L)$ be given by the formulas:
$$T(f)\restriction_{2^{\omega}\cup\{\infty\}}=f\restriction_{2^{\omega}\cup\{\infty\}},$$
and
$$T(f)(n)=\frac{f(\sigma^{-1}(n))}{n+1}$$
for $n\in\omega$, and
$T(f)(\omega)=0$. Since $f\in C_{p}(K)$ where $K$ is compact, the function $f$ is bounded, and therefore $T(f)(n)\xrightarrow[]{}0$ when $n\xrightarrow[]{}\infty$, so $T(f)\in C_{p}(L)$ for every $f\in C_{p}(K)$. It is easy to see that $T$ is continuous, linear, injective and open onto its image, and hence the transformation $T$ is a homeomorphic embedding.

\end{proof}
\section{A characterisation of $\kappa$-Corson compact spaces in terms of function spaces $C_{p}(K)$}
\label{char kappa-cors sekcja}
In this section we present a generalisation of results due to Pol \cite{P} and Bell and Marciszewski \cite{BM}. Pol's theorem characterises Corson compacta $K$ in terms of function spaces $C_{p}(K)$. Bell and Marciszewski extended this result to $\kappa^{+}$-Corson compact spaces. We also conclude, using results from \cite{MPZ} 
 and Sections \ref{secdimension}, \ref{secNY} that the class of $\omega$-Corson compact spaces $K$ is invariant under linear homeomorphisms of function spaces $C_{p}(K)$.

Given an infinite cardinal numbers $\kappa$ and $\lambda$, let $L_{\kappa}(\lambda)$
denote the set $\lambda\cup\{\infty\}$ topologized as follows: all points  $\alpha\in\lambda$ are isolated, and open neighbourhoods 
of $\infty$ are of the form $\{\infty\}\cup A$ where $A\subseteq \lambda$ and $|\lambda\setminus A|<\kappa$.  
For an infinite cardinal number $\kappa$, let $\mathscr{L}_{\kappa}$ denote the class of all spaces which are continuous images of closed subsets of the countable product $L_{\kappa}(\lambda)^{\omega}$. Pol's result states that a compact space $K$ is Corson compact if and only if $C_{p}(K)\in\mathscr{L}_{\omega_{1}}$. Bell and Marciszewski extended this theorem to $\kappa^{+}$-Corson compact spaces. 

We generalise the results from \cite{BM} and \cite{P} to $\kappa$-Corson compact spaces for any  regular, uncountable cardinal number $\kappa$. Our reasoning is a modification of the proof from \cite{BM} which in turn follows the arguments from \cite{P}.

Let $\kappa$ be a regular, uncountable cardinal number, and let $\lambda\in Card$, $\lambda>\kappa$. For $A\subseteq L_{\kappa}(\lambda)$ such that $\infty\in A$, define the retraction $r_{A}:L_{\kappa}(\lambda)\xrightarrow[]{}A$ 
$$r_{A}(x)=\begin{cases}
x & \text{if } x\in A\\
\infty & \text{if } x\notin A  
\end{cases}.$$
Let $R_{A}=\prod_{n\in\mathbb{N}} r_{A}:L_{\kappa}(\lambda)^{\omega}\xrightarrow[]{}A^{\omega}$. For $n\in\mathbb{N}$, let $P_{n}:L_{\kappa}(\lambda)^{\omega}\xrightarrow[]{}L_{\kappa}(\lambda)^{\omega}$ be given by the formula $P_{n}(x_{1},x_{2},\dots)=(x_{1},\dots x_{n},\infty,\infty,\dots).$

\begin{Lemma}
Let $F\subseteq L_{\kappa}(\lambda)^{\omega}$ be a closed subset, and let $A$ be a subset of $L_{\kappa}(\lambda)$, such that $\infty\in A$ and $|A|\geq \kappa$. Then there exists $B\subseteq L_{\kappa}(\lambda)$ such that $A\subseteq B$, $|B|=|A|$ and $R_{B}(F)\subseteq F$. \label{powiększanie dużych zbiorów}
\end{Lemma}
\begin{proof}
For every $z\in L_{\kappa}(\lambda)^{\omega}\setminus F$, fix a basic open neighbourhood 
$$V^{z}=V_{1}^{z}\times V_{2}^{z}\times\dots V_{k}^{z}\times L_{\kappa}(\lambda)\times\dots$$
of $z$, disjoint from $F$ such that $k=k(z)$ is minimal for all such neighbourhoods. Define 
$$C_{z}=\bigcup\{L_{\kappa}(\lambda)\setminus V_{i}^{z}:\infty\in V_{i}^{z},i\leq k\},$$ and observe that $|C_{z}|<\kappa$. We construct the set $B$ by induction. Start with $B_{1}=A$, and put 
$$B_{n+1}=\bigcup\{C_{z}:z\in P_{n}(B_{n}^{\omega})\setminus F\}\cup B_{n},\:n\in\mathbb{N}. $$ 
Finally, define $B=\bigcup_{n\in\mathbb{N}}B_{n}$. By induction, if $|B_{n}|=|A|$, then $|P_{n}(B_{n}^{\omega})|=|A|$, so $|C_{z}|<\kappa$ implies $|B_{n+1}|=|A|$. Then, we have $|B|=|A|$. We will show that $R_{B}(F)\subseteq F$.

Let $x=(x_{i})\in F$, and suppose that $y=(y_{i})=R_{B}(x)\notin F$. Notice that $y_{i}\in B$ for every $i\in\mathbb{N}$. Take $z=(z_{i})=P_{k}(y)$, where $k=k(y)$. Then $z_{i}=y_{i}$ for $i\leq k$, hence $z\notin F$ and $k(z)=k$. Choose $n\geq k$ such that $z_{i}\in B_{n}$, for $i\leq k$. Observe that $C_{z}\subseteq B$. For each $i\leq k$, we have $z_{i}=y_{i}=r_{B}(x_{i})$, and we can consider two possibilities: Firstly, if $z_{i}\neq\infty$, then $x_{i}=r_{B}(x_{i})$ by the definition of $r_{B}$, so $x_{i}=z_{i}\in V_{i}^{z}$. Secondly, if $z_{i}=\infty$, then $x_{i}\notin C_{z}$. Indeed, if $x_{i}\in C_{z}\subseteq B$, then $z_{i}=r_{B}(x_{i})=x_{i}\neq\infty$, since $\infty\notin C_{z}$. Therefore again $x_{i}\in V_{i}^{z}$. This shows that $x$ belongs to $V^{z}$ which is disjoint from F, contradiction.   
\end{proof}
\begin{Lemma}\label{powiększanie małych zbiorów}
Let $F\subseteq L_{\kappa}(\lambda)^{\omega}$ be a closed subset, and let $A$ be a subset of $L_{\kappa}(\lambda)$ such that $\infty\in A$ and $|A|< \kappa$. Then there exists $B\subseteq L_{\kappa}(\lambda)$ such that $A\subseteq B$, $|B|<\kappa$ and $R_{B}(F)\subseteq F$.
\end{Lemma}
\begin{proof}
Let us define the set $B$ the same way as in the proof of Lemma \ref{powiększanie dużych zbiorów}. We will show that $|B|<\kappa$. Firstly, we have $|B_{2}|<\kappa$ as it is a union of less than $\kappa$ sets of cardinality less than $\kappa$. It is easy to see that, by induction, $|B_{n}|<\kappa$ for every $n\geq 2$. Finally, we have $|B|=|\bigcup_{n\in\mathbb{N}}B_{n}|<\kappa$ because $\kappa$ is regular. The inclusion $R_{B}(F)\subseteq F$ can be shown in the same way as in Lemma \ref{powiększanie dużych zbiorów}.
\end{proof}
\begin{Lemma}
\label{Ciąg zbiorów}
Let $F\subseteq L_{\kappa}(\lambda)^{\omega}$ be a closed subset, where $\lambda\geq \kappa$. There exists a family $\{A_{\alpha}:\alpha\leq\lambda\}$ of subsets of $L_{\kappa}(\lambda)$ such that:
\begin{enumerate}[(1)]
\item $\infty\in A_{0}$,
\item $|A_{\alpha}|<\kappa$, for $\alpha<\kappa$,
\item $|A_{\alpha}|=|\alpha|$, for $\kappa\leq\alpha\leq\lambda$,
\item $A_{\alpha}\subseteq A_{\beta}$, for $\alpha\leq\beta\leq\lambda,$
\item $A_{\alpha}=\bigcup_{\beta<\alpha}A_{\beta}$, for $\alpha\in \Limm\cap(\lambda+1),$
\item $A_{\lambda}=L_{\kappa}(\lambda),$
\item $R_{A_{\alpha}}(F)\subseteq F$ for $\alpha\leq\lambda.$
\end{enumerate}
\end{Lemma}
\begin{proof}
By Lemma \ref{powiększanie małych zbiorów} for $F$ and $A=\{\infty\}$, there exists $A_{0}\subseteq L_{\kappa}(\lambda)$ such that $\infty\in A_{0}$,\\ $|A_{0}|<\kappa$ and $R_{A_{0}}(F)\subseteq F$. We will construct sets $A_{\alpha}$ for $\alpha\geq 1$ by induction. Assume we have constructed $A_{\alpha}$ for $\alpha\leq\beta<\kappa$. By Lemma \ref{powiększanie małych zbiorów}, there exists $B\subseteq L_{\kappa}(\lambda)$ such that $$A_{\beta}\cup(\beta+1)\subseteq B, R_{B}(F)\subseteq F \text{ and } |B|<\kappa.$$ Let $A_{\beta+1}=B$. For $\beta\in Lim\cap(\kappa+1)$, put 
$A_{\beta}=\bigcup_{\alpha<\beta}A_{\alpha}$. As  
$\beta+1\subseteq A_{\beta+1}$ for $\beta<\kappa$, we have 
$|A_{\kappa}|=\kappa$.

Assume that we constructed $A_{\alpha}$ for 
$\alpha\leq\beta<\lambda$ where $\beta\geq\kappa$. By Lemma \ref{powiększanie dużych zbiorów}, there exists $B\subseteq L_{\kappa}(\lambda)$, such that $R_{B}(F)\subseteq F$, $|B|=|A_{\beta}|$ and $A_{\beta}\cup (\beta+1)\subseteq B$. For a limit ordinal $\kappa<\beta\leq\lambda$, let $A_{\beta}=\bigcup_{\alpha<\beta}A_{\alpha}$.

It remains, to check that if $R_{A_{\alpha}}(F)\subseteq F$ for $\alpha<\beta$ where $\beta\in Lim$, then $R_{A_{\beta}}(F)\subseteq F$. Assume that there exists $x\in F$ such that $R_{A_{\beta}}(x)\notin F$. Let  $V=\prod_{i\in\mathbb{N}} V_{i}$ a basic neighbourhood of $R_{A_{\beta}}(x)$ such that $V_{i}=L_{\kappa}(\lambda)$ for $i> n$ and $V\cap F=\emptyset$. For every \mbox{$k\in\{1,2,\dots,n\}$,} there exists $\alpha_{k}\leq\beta$ such that for every $\alpha_{k}<\alpha<\beta$, we have \mbox{$r_{A_{\alpha_{k}}}(x_{k})=r_{A_{\alpha}}(x_{k})$.} Consequently, the point $R_{A_{\alpha}}(x)$ does not belong to $F$ for $\alpha$ greater than $max(\alpha_{1},\dots,\alpha_{n})$, contradiction.
\end{proof}
Let $\alpha$ be a limit ordinal number. Observe that for every $x\in A_{\alpha}$ there exists a minimal $\alpha_{0}<\alpha$ for which $x\in A_{\alpha_{0}}$, so for $\alpha_{0}\leq\beta\leq\alpha$, we have $r_{\alpha_{0}}(x)=r_{\beta}(x)$. This shows that, for a limit ordinal $\alpha$ and $y\in L_{\kappa}(\lambda)^{\omega}$,
\begin{equation}
\label{granica}
R_{A_{\alpha}}(y)=\Lim{\beta\xrightarrow[]{}\alpha}R_{A_{\beta}}(y).
\end{equation}
  Let $\kappa$ be an infinite cardinal number. Recall that a space is called $\kappa$-Lindel\"of if every open cover has a subcover of cardinality less than $\kappa$. 
The following proposition follows from \\ \cite[Corollary 4.2]{No}.
\begin{prop}\label{liczba Lindelofa} 
Let $\kappa$ be an infinite regular cardinal number. Every space from the class $\mathscr{L}_{\kappa}$ is $\kappa$-Lindel\"of.
\end{prop}
\begin{prop}\label{domkn podzb ci wloz CP} 
For a closed set $F\subseteq L_{\kappa}(\lambda)^{\omega}$, there exists a continuous, linear injection \\$T:C_{p}(F)\xrightarrow[]{} 
\Sigma_{\kappa}(\mathbb{R}^{\lambda})$.
\end{prop}
\begin{proof}
If $\lambda<\kappa$, then $L_{\kappa}(\lambda)$ is discrete, and so $F$ is metrizable, and $w(F)\leq \lambda$. There exists a dense $D\subseteq F$ such that $|D|=\lambda$. Let $T:C_{p}(F)\xrightarrow[]{} \mathbb{R}^{D}$ be the restriction map $T(f)=f\restriction_{D}$. Clearly $T$ is linear, injective and continuous, and $R^{D}$ embeds in $\Sigma_{\kappa}(\mathbb{R}^{\lambda})$ as $|D|<\kappa$.

Assume $\lambda\geq\kappa$, and assume that the thesis holds for every $\mu<\lambda$. Let $F$  be a closed subset of $L_{\kappa}(\lambda)^{\omega}$, and let the family $\{A_{\alpha}:\alpha\leq\lambda\}$ be as in Lemma \ref{Ciąg zbiorów}. Let 
$R_{\alpha}=R_{A_{\alpha}}\restriction_{F}$ and $F_{\alpha}=R_{\alpha}(F)$. 
Function $G:C_{p}(F_{\alpha})\xrightarrow[]{}C_{\alpha}$ given by the formula $G(f)=f\circ R_{\alpha}$ where 
$$C_{\alpha}=\{f\circ R_{\alpha}:f\in C_{p}(F_{\alpha})\}\subseteq C_{p}(F)$$
is a linear homeomorphism. 
Since $A_{\alpha}\subseteq A_{\beta}$ for $\alpha<\beta$, we have $R_{\alpha}=R_{\alpha}\circ 
R_{\beta}$. This implies $C_{\alpha}\subseteq C_{\beta}$ for $\alpha<\beta$.

We will check that 
$F_{\alpha}=F\cap A_{\alpha}^{\omega}$. Firstly, $F_{\alpha}=R_{\alpha}(F)\subseteq 
A_{\alpha}^{\omega}$ as $im\:(r_{A_{\alpha}})\subseteq A_{\alpha}$ because 
$\infty\in A_{\alpha}$. Secondly, we have $R_{A_{\alpha}}(F)\subseteq F$ (condition (7) from Lemma \ref{Ciąg zbiorów}). Finally, if $x\in F\cap A_{\alpha}^{\omega}$, then 
$R_{A_{\alpha}}(x)=x$, so $R_{A_{\alpha}}(x)=x\in F$, therefore $x\in
R_{A_{\alpha}}(F)=F_{\alpha}$.

Since $|A_{\alpha}|<\lambda$, $F_{\alpha}\subseteq L_{\kappa}(A_{\alpha})^{\omega}$ and  $C_{p}(F_{\alpha})\sim C_{\alpha}$ for $\alpha<\lambda$, by inductive hypothesis, there exists a continuous linear injection 
$T_{\alpha}:C_{\alpha}\xrightarrow[]{}
\Sigma_{\kappa}(\mathbb{R}^{|A_{\alpha}|})$.  For the special case $\lambda=\kappa$, define $T_{\kappa}:C_{p}(F_{\kappa})\xrightarrow[]{}\prod_{\alpha\in\kappa}\Sigma_{\kappa}(\mathbb{R}^{|A_{\alpha}|})$ as
$$T_{\kappa}(f)=(T_{0}(f\circ R_{0}),T_{\alpha+1}(f\circ R_{\alpha+1}-f\circ R_{\alpha})_{\alpha\in[0,\kappa)}),$$ 
and $T_{\lambda}:C_{p}(F_{\lambda})\xrightarrow[]{}\prod_{\alpha\in\lambda}\Sigma_{\kappa}(\mathbb{R}^{|A_{\alpha}|})$ for $\lambda>\kappa$ as
$$T_{\lambda}(f)=(T_{\kappa}(f\circ R_{\kappa}),(T_{\alpha+1}(f\circ R_{\alpha+1}-f\circ R_{\alpha}))_{\alpha\in[\kappa,\lambda)}).$$
It is clear that $T_{\kappa}$ is continuous and linear. Let us check that it is injective. Assume that $T_{\kappa}(f)=\mathbf{0}$ for some nonzero $f\in C_{p}(F_{\kappa})$. Let $\alpha=min\{\beta\leq\kappa:f\circ R_{\beta}\neq\mathbf{0}\}$ (such $\beta$ exists since $R_{\kappa}=id_{F}$).  If $\alpha=0$, then $T_{0}(f\circ R_{0})\neq\mathbf{0}$ because $T_{0}$ is injective, so we have a contradiction. If $\alpha$ is a successor ordinal, then let $\alpha=\beta+1$. In this case
$$T_{\beta+1}(f\circ R_{\beta+1}-f\circ R_{\beta})=T_{\beta+1}(f\circ R_{\beta+1})-T_{\beta+1}(f\circ R_{\beta})\neq \mathbf{0}$$
because $T_{\beta+1}(f\circ R_{\beta})=\mathbf{0}$ by the definition of $\alpha$, and $T_{\beta+1}(f\circ R_{\beta+1})\neq\mathbf{0}$ by injectivity of $T_{\beta+1}$.
If $\alpha$ is a limit ordinal number, then by Equation \ref{granica}, we arrive at contradiction. 

We can identify $\prod_{\alpha\in\kappa}\Sigma_{\kappa}(\mathbb{R}^{|A_{\alpha}|})$ with a subspace of $\mathbb{R}^{S}$ where $S=\dot\bigcup_{\alpha\in\kappa} A_{\alpha}$, $|S|=\kappa$.
Now we will check that for every $f\in C_{p}(F)$, $|supp(T(f))|<\kappa$. By regularity of $\kappa$, it is enough to show that 
$$T_{\alpha+1}(f\circ R_{\alpha+1}-f\circ R_{\alpha})\neq \mathbf{0}$$
for less than $\kappa$ ordinal numbers $\alpha$ because $|A_{\alpha}|<\kappa$ for $\alpha<\kappa$. Let $f\in C_{p}(F)$ and suppose to
the contrary that we can find a set $T\subseteq S$ of cardinality $\kappa$ such that, for every $\alpha\in T$, there exists 
$x^{\alpha}=(x^{\alpha}_{i})\in F$ for which 
$f(R_{\alpha+1}(x^{\alpha}))-f(R_{\alpha}(x^{\alpha}))\neq \mathbf{0}.$
Since $cf(\kappa)=\kappa>\omega$, we can assume that, for 
every $\alpha\in T$, we have $$|f(R_{\alpha+1}(x^{\alpha}))-f(R_{\alpha}(x^{\alpha}))|>\epsilon$$ for some $\epsilon>0$. Using Proposition \ref{liczba Lindelofa} for $F$, we can find a 
point \mbox{$x=(x_{i})\in F$} such that for every neighbourhood $U$ of $x$ we have $|{\alpha\in T:R_{\alpha}(x^{\alpha})\in U}|=\kappa$. Let 
$$V=V_{1}\times\dots\times V_{k}\times L_{\kappa}(\lambda)\times\dots$$ be a basic neighbourhood of $x$ in $L_{\kappa}(\lambda)^{\omega}$ such that $diam (f(V\cap F))<\epsilon$ 
and $V_{i}=\{x_{i}\}$ if $x_{i}\neq\infty$.

Let 
$$C=\bigcup\{L_{\kappa}(\lambda)\setminus V_{i}:\infty\in V_{i},i\leq k\}.$$
Since $|C|<\kappa$, and we have $\kappa$ 
many points $R_{\alpha}(x^{\alpha})$ in $V$ with $\alpha\in T$, we can find $\alpha\in T$ such that $R_{\alpha}(x^{\alpha})\in V$ and the set $A_{\alpha+1}\setminus A_{\alpha}$ is
disjoint from $C$. Then we also have $R_{\alpha+1}(x^{\alpha})\in V$. Indeed, we can show that $R_{\alpha+1}(x^{\alpha})_{i}\in V_{i}$ for all $i\leq k$.

If $\infty\notin V_{i}$, then $R_{\alpha}(x^{\alpha})_{i}=r_{A_{\alpha}}(x_{i}^{\alpha})\in V_{i}$ and $\infty\notin V_{i}$, so $r_{A_{\alpha}}(x_{i}^{\alpha})=x_{i}=x_{i}^{\alpha}\in A_{\alpha}\subseteq A_{\alpha+1}$. Hence $R_{\alpha+1}(x^{\alpha})_{i}=r_{A_{\alpha+1}}(x^{\alpha}_{i})=x_{i}^{\alpha}=x_{i}\in V_{i}$.

If $\infty\in V_{i}$, then $(A_{\alpha+1}\setminus A_{\alpha})\cap C=\emptyset$ implies that $A_{\alpha+1}\setminus A_{\alpha}\subseteq V_{i}$. There are two cases. Firstly, if $x_{i}^{\alpha}\in A_{\alpha+1}\setminus A_{\alpha}$, then $R_{\alpha+1}(x^{\alpha})_{i}=r_{A_{\alpha+1}}(x_{i}^{\alpha})=x_{i}^{\alpha}\in V_{i}$. Secondly, if $x_{i}^{\alpha}\notin A_{\alpha+1}\setminus A_{\alpha}$, then $R_{\alpha+1}(x^{\alpha})_{i}=R_{\alpha}(x^{\alpha})_{i}=r_{A_{\alpha}}(x_{i}^{\alpha})\in V_{i}$.

By the definition of $V$, we have a contradiction with $|f(R_{\alpha+1}(x^{\alpha}))-f(R_{\alpha}(x^{\alpha}))|>\epsilon$. Analogous reasoning yields that $T_{\lambda}$ is injective and $im\:T_{\lambda}\subseteq\Sigma_{\kappa}(\mathbb{R}^{\lambda})$ for $\lambda>\kappa$.
\end{proof}
\begin{prop}
If $C_{p}(K)\in\mathscr{L}_{\kappa}$, then $K$ is a $\kappa$-Corson compact space.
\label{Cp implikuje kappa corson}
\end{prop}
\begin{proof}
Let $F\subseteq L_{\kappa}(\lambda)^{\omega}$ be a closed subset. Assume that  there is a continuous surjection $\phi:F\xrightarrow[]{}C_{p}(K)$. Define $\phi':C_{p}(C_{p}(K))\xrightarrow[]{}C_{p}(F)$ as $\phi'(g)=g\circ\phi$. It is a well known fact, that in this situation $\phi'$ is a continuous injection. Let $\psi:K\xrightarrow[]{}C_{p}(C_{p}(K))$ be the standard embedding given by evaluations.
By Proposition \ref{domkn podzb ci wloz CP}, there exists a continuous
linear injection $\theta:C_{p}(F)\xrightarrow[]{}
\Sigma_{\kappa}(\mathbb{R}^{\Gamma})$ for some $\Gamma$. The function $\theta\circ\phi'\circ \psi:K\xrightarrow[]{}\Sigma_{\kappa}(\mathbb{R}^{\Gamma})$ is a continuous injection from a compact, Hasudorff space into a Hausdorff space, therefore it is a homeomorphic embedding. This shows that $K$ is a $\kappa$-Corson compact space.
\end{proof}
The theorem below can be proved using the same arguments as in the proof of Theorem 6.1 in \cite{BM}.
\begin{Theorem}
Let $\kappa$ be an uncountable cardinal number. If space $K$ is $\kappa$-Corson compact, then $C_{p}(K)\in\mathscr{L}_{\kappa}$.
\label{kappa corson implikuje Cp}
\end{Theorem}
By Proposition \ref{Cp implikuje kappa corson} and Theorem \ref{kappa corson implikuje Cp}, we get
\begin{Theorem}
For a compact space $K$ and a regular $\kappa>\omega$, the space $K$ is $\kappa$-Corson compact if and only if $C_{p}(K)\in\mathscr{L}_{\kappa}$.
\label{char kappa - Corson}
\end{Theorem}
\begin{Theorem}
Let $\kappa$ be a regular, uncountable cardinal number. Assume that $K$ is a \mbox{$\kappa$-Corson} compact space and $L$ is compact. If there exists a continuous surjective transformation $T:C_{p}(K)\xrightarrow[]{}C_{p}(L)$, then $L$ is $\kappa$-Corson compact as well.
\end{Theorem}
\begin{proof}
By Theorem \ref{char kappa - Corson}, we have $C_{p}(K)\in\mathscr{L}_{\kappa}$. There exists $\lambda\in\,Card$, a closed subset $F$ of $ L_{\kappa}(\lambda)^{\omega}$ and a continuous surjection $f:F\xrightarrow[]{}C_{p}(K)$. The composition $T \circ f:F\xrightarrow[]{}C_{p}(L)$ is a continuous surjection as well, and therefore by Theorem \ref{char kappa - Corson}, the space $L$ is $\kappa$-Corson compact.
\end{proof}
The theorem below was originally proved in \cite[Corollary 2.17]{BKT}
\begin{Theorem}
For a regular, uncountable cardinal number $\kappa$ the class of $\kappa$-Corson compact spaces is invariant under continuous images.
\label{stable}
\end{Theorem}
\begin{proof}
Let $K$ be a $\kappa$-Corson compact space for some regular, uncountable $\kappa$. By Theorem \ref{char kappa - Corson}, we have $C_{p}(K)\in\mathscr{L}_{\kappa}$. If $L$ is a continuous image of $K$, then there exists a closed homeomorphic embedding of $C_{p}(L)$ into $C_{p}(K)$, consequently $C_{p}(L)\in\mathscr{L}_{\kappa}$. By Theorem \ref{char kappa - Corson}, the space $L$ is $\kappa$-Corson compact.
\end{proof}
The theorem below is \cite[Corollary 5.1]{MPZ}
\begin{Theorem}
    An $NY$ compact space $K$ is \mbox{$\omega$-Corson} compact if and only if it is strongly countable dimensional.   \label{NY w-corson sil przel}
\end{Theorem}
By Theorems \ref{sil kappa wymiar zach przy lin homeo}, \ref{NY zachowuja sie przy lin homeo},  \ref{NY w-corson sil przel}, we obtain
\begin{Theorem}
    Let $K$ and $L$ be compact spaces. Assume there exists a continuous, linear transformation $\phi:C_{p}(K)\xrightarrow[]{}C_{p}(L)$ such that $\phi(C_{p}(K))$ is dense in $C_{p}(L)$. If $K$ is \mbox{$\omega$-Corson} compact, then $L$ is $\omega$-Corson compact as well.
    \label{omega corson zach przy lin homeo}
\end{Theorem}

\end{document}